# Introducing an Analysis in Finite Fields

Hélio M. de Oliveira and Ricardo M. Campello de Souza

Universidade Federal de Pernambuco UFPE
C.P.7800 - 50.711-970 - Recife - PE, Brazil
E-mail: {hmo,ricardo}@ufpe.br

*Abstract - Looking forward to introducing an analysis in Galois Fields, discrete functions are considered (such as transcendental ones) and MacLaurin series are derived by Lagrange's Interpolation. A new derivative over finite fields is defined which is based on the Hasse Derivative and is referred to as negacyclic Hasse derivative. Finite field Taylor series and $\alpha$-adic expansions over GF(p), p prime, are then considered. Applications to exponential and trigonometric functions are presented. Theses tools can be useful in areas such as coding theory and digital signal processing.*

## 1. INTRODUCTION

One of the most powerful theories in mathematics is the Real Analysis [1]. This paper introduces new tools for Finite Fields which are somewhat analogous to classical Analysis. Finite fields have long been used in Electrical Engineering [2], yielding challenging and exciting practical applications including deep space communication and compact discs (CD). Such a theory is the basis of areas like algebraic coding, and is related, to various degrees, to topics such as cryptography, spread-spectrum and signal analysis by finite field transforms [3,4]. Some new and unexpected applications of this framework are currently under developing which include bandwidth efficient code division multiplex [5]. Hereafter we adopt the symbol := to denote *equal by definition*, the symbol $\equiv$ to denote a *congruence mod p* and K[$x$] a *polynomial ring over a field* K.

## 2. FUNCTIONS OVER GALOIS FIELDS

For the sake of simplicity, we consider initially Galois Fields GF($p$), where $p$ is an odd prime.

*Functions:* Let $\alpha$ be an element of GF($p$) of order $N$.
i) One of the simplest mapping is the "affine function" which corresponds to a straight line on a finite field, e.g. a($x$)=3$x$+4, $x \in$ GF(5).
$x$=  0  1  2  3  4
a($x$)=  4  2  0  3  1.
There exists $p(p-1)$ linear GF($p$)-functions.
ii) Exponential function: $\alpha^i$, $i$=0,1,...,$N$-1.
e.g. letting $\alpha$=2, the exponential function $2^x$ corresponds to $2^i$ (mod 5).
$i$=  0  1  2  3  4       indexes $i$ mod 5
       1  2  4  3  1.
It is also interesting to consider "*shortened*" functions in which indexes are reduced modulo $p-k$ for 0<$k$<$p$.

$i$=  0  1  2  3  |  4       indexes $i$ mod 4
       1  2  4  3       repeat.
Notice that $(\forall x)\ 2^x \neq 0$ and $(\forall i)\ 2^i \not\equiv 0$.

Composite functions can also be generated as usual.
For instance, the function
$i$=  0  1  2  3  4
       1  4  1  3  2
corresponds to $2^{3x+4}$ (mod 5) = $2^x \circ (3x+4)$.

Therefore inverse functions over GF($p$) can be defined as usual, i.e., a function a($x$) has an inverse a$^{-1}$($x$) iff a($x$) o a$^{-1}$($x$)$\equiv x$ mod $p$.
There are $p!$ inversive GF($p$)-functions. The inverse of a($x$)$\equiv$3$x$+4 is a$^{-1}$($x$)$\equiv$2$x$+2 mod 5.
Since the characteristic of the field is different from 2, we can define the odd and even component of a function, respectively, as:

$$o(x) := \frac{a(x) - a(-x)}{2} \text{ and } e(x) := \frac{a(x) + a(-x)}{2}.$$

An even (respectively odd) function has $o(x) \equiv 0$ (respectively $e(x) \equiv 0$).
The $2^x$ exponential function defined over GF(5) has
$x$=  0  1  2  3  4    $x$=  0  1  2  3  4
$o(x)$=  3  3  3  2  2.  $e(x)$=  3  4  1  1  4.

iii) $k$-Trigonometric functions
A trigonometry in Finite Fields was recently introduced [6]. Let G($q$) be the set of Gaussian integers over GF($q$),i.e. G($q$):={$a+jb$, $a,b \in$ GF($q$)}, where $q=p^r$, $p$ being an odd prime, and $j^2$=-1 is a quadratic non-residue in GF($q$). With the usual complex number operations, we have then that GI:= <G($q$),$\oplus$,$\otimes$,> is a field.

Definition 1.
Let $\alpha$ have multiplicative order $N$ in GF($q$). The GI-valued $k$-trigonometric functions cos$_k$(.) and sin$_k$(.) are:

$$\cos_k(\alpha^i) := \frac{1}{2}(\alpha^{ik} + \alpha^{-ik}),$$

$$\sin_k(\alpha^i) := \frac{1}{2j}(\alpha^{ik} - \alpha^{-ik}) \quad i,k=0,1,...N-1.$$

The properties of the unit circle hold, i.e.,
$$\sin_k^2(i) + \cos_k^2(i) \equiv 1 \mod p.$$

We consider now $k$-trigonometric functions with $\alpha$=3 an element of order 6 in GF(7). A pictorial representation of the cos-function on this finite field is presented in the sequel. There exists several $k$-cos: {cos$_k$($i$)}, $k$=0,1,..,5. In order to clarify the symmetries, half of the elements of GF($p$) are considered positive and the other half assumes negative values. Therefore, we represent the Galois field elements as the set $\{0, \pm 1, \pm 2, ..., \pm(p-1)/2\}$.

The $\cos_k(i)$ over GF(7) assumes the following values ($k,i=0,...,5$):

$\cos_0(i)$=1 1 1 1 1 1 ≡ 1 1 1 1 1 1 (mod 7)

$\cos_1(i)$=1 4 3 6 3 4 ≡ 1 -3 3 -1 3 -3 (mod 7)

$\cos_2(i)$=1 3 3 1 3 3 ≡ 1 3 3 1 3 3 (mod 7)

$\cos_3(i)$=1 6 1 6 1 6 ≡ 1 -1 1 -1 1 -1 (mod 7)

$\cos_4(i)$=1 3 3 1 3 3 ≡ $\cos_2(i)$

$\cos_5(i)$=1 4 3 6 3 4 ≡ $\cos_1(i)$

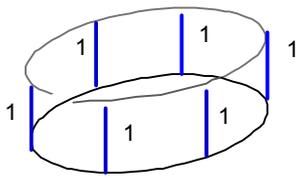

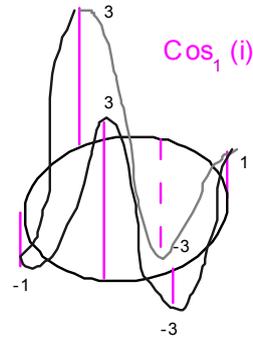

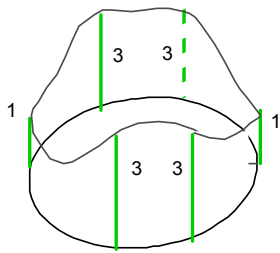

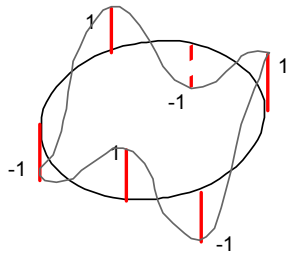

Figure 1. Pictorial representation of $\cos_k(i)$ in GF(7).

The "dough-nut ring" representation is well suitable for handling with the cyclic structure of finite fields. The envelope does not exists: It just gives some insight on the concept of different "frequencies" on each carrier (oscillatory behavior). We define $\pi$ on GF($p$) as

$(\pi):=(p-1)/2$. Clearly, $\cos_1(i+(\pi))=-\cos_1 i$.

We shall see that any function from GF($p$) to GF($p$) defined by $N \mid p$-1 points can be written as a polynomial of degree $N$-1 (i.e., $N$ indeterminates), e.g., $2^i \equiv i^3 + 1$ (mod 5) (indexes $i$ mod 4).
Generally there are $N$ unknown coefficients, $f(i) \equiv a_0+a_1i+a_2i^2+...+a_{N-1}i^{N-1}$ mod $p$, and such a decomposition corresponds to the (finite) MacLaurin series with the advantage that there are no errors.

## 3. DETERMINING THE SERIES

Given a function $f(i)$, the only polynomial "passing through" all pair of points $\{i, f(i)\}$ can be found by solving a Linear System over GF($p$) [7]. For instance, the 1-cos function over GF(5) leads to:

$$\begin{bmatrix} 1 & 0 & 0 & 0 \\ 1 & 1 & 1 & 1 \\ 1 & 2 & 4 & 3 \\ 1 & 3 & 4 & 2 \end{bmatrix} \begin{pmatrix} a_0 \\ a_1 \\ a_2 \\ a_3 \end{pmatrix} \equiv \begin{bmatrix} 1 \\ 0 \\ 4 \\ 0 \end{bmatrix} \text{ (mod 5)}.$$

Another way to find out the coefficients of the polynomial expansion is by using Lagrange's interpolation formula [8],

$$L_3(i) \equiv 1 \frac{(i-1)(i-2)(i-3)}{(0-1)(0-2)(0-3)} +$$
$$4 \frac{(i-0)(i-1)(i-3)}{(2-0)(2-1)(2-3)} \quad \text{(mod 5)}.$$

Therefore
$$\cos_1(i) \equiv 1 - 2i - i^2 + 2i^3 \text{ (mod 5)}.$$

Identically, the 1-sin function [6] can be expanded:
$i$=      0   1   2   3
$\sin_1(i)$    0   $j3$   0   $j2$,

$$\sin_1(i) \equiv j(i^3 - i^2 - 2i) \quad \text{(mod 5)}.$$

Euler's formula over a Galois field proved in [6] can also be verified in terms of series:
$\cos(i)+j\sin(i) \equiv$
$(1-2i-i^2+2i^3)-(i^3-i^2-2i) \equiv i^3+1 \equiv 2^i$ mod 5 ($i$ mod 4).
Those results hold despite $\sqrt{-1} \equiv 2$ mod 5!
The unicity of the series decomposition can be established by:

Proposition 1. Given a function f defined by its values $f(x)$, $\forall x \in$ GF($p$), there exist only one Maclaurin series for $f$.

Proof. Letting $f(x) \equiv a_0+a_1x+a_2x^2+...+a_{p-1}x^{p-1}$ (mod $p$) then $f(0)=a_0$ and we have the following linear system

$$\begin{bmatrix} 1 & 1 & 1 & \ldots & 1 \\ 2 & 2^2 & 2^3 & \ldots & 2^{p-1} \\ 3 & 3^2 & 3^3 & \ldots & 3^{p-1} \\ \vdots & \ldots & \ldots & \ldots & \vdots \\ (p-1) & (p-1)^2 & (p-1)^3 & \ldots & (p-1)^{p-1} \end{bmatrix}.$$

$$\begin{pmatrix} a_1 \\ a_2 \\ a_3 \\ \vdots \\ a_{p-1} \end{pmatrix} \equiv \begin{pmatrix} f(1)-f(0) \\ f(2)-f(0) \\ f(3)-f(0) \\ \vdots \\ f(p-1)-f(0) \end{pmatrix} \pmod{p}$$

or $V \cdot \underline{a} = \underline{f}$. Since V is a Vandermonde matrix, its determinant is:

$\det[V] \equiv \prod_{i<j}(x_j - x_i) \pmod{p}$ where $x_i := i$.

Therefore, $\det[V] \not\equiv 0 \mod p$ concluding the proof. ∎

## 4. DERIVATIVES ON FINITE FIELDS

Over a finite field, we are concerned essentially with derivative of polynomials since other functions can be expanded in a series. We denote by $a^{(i)}(x)$ the $i$-th derivative of the polynomial $a(x)$. The concept of classical derivative however presents a serious drawback, namely the derivatives of order greater or equal to the characteristic of the field vanish. Let $a(x) = \sum_{i=0}^{N-1} a_i x^i$ where $(\forall i)\ a_i \in GF(p)$, $x$ is a dummy variable. Then the derivatives $a^{(p)}(x) \equiv 0$ and $a^{(i)}(x) \equiv 0$ ($\forall i \geq p$) no matter the coefficients. A more powerful concept of derivative over a finite field was introduced a long time ago (see [9]):

<u>Definition 2</u>. The Hasse derivative of a polynomial function $a(x) = \sum_{i=0}^{N-1} a_i x^i$ is defined by

$a^{[r]}(x) := \sum_{i=0}^{N-1} \binom{i}{r} a_i x^{i-r}$ with $\binom{i}{r} \equiv 0$ for $i<r$.

Clearly the classical Newton-Leibnitz derivative yields

$a^{(r)}(x) = \sum_{i=0}^{N-1} i(i-1)\ldots(i-r+1) a_i x^{i-r}$. It follows that $\frac{1}{r!} a^{(r)}(x) = a^{[r]}(x)$.

The Hasse derivative has been successfully used on areas where finite fields play a major role such as coding theory [10].

In the following, we introduce a new concept of derivative in a finite field. In order to exploit its cyclic structure, a polynomial ring $GF(p)[x]$ is required. Let 1 be the constant function $f(x)=1\ \forall x \in GF(p)$.

Conventional Hasse derivative considers polynomial modulo $x^p - 1$ so that $\frac{d1}{dx} \equiv 0$. On the other hand, the fact that $\binom{i}{r} = 0\ \forall i<r$ does not allow to handle with negative powers of $x$.

In contrast with Hasse derivative which always decreases the degree of the polynomial, we introduce a new concept of negacyclic Hasse derivative considering the polynomial ring $GF(p)[x]$ reduced modulo $x^{p-1}+1$ [11]. Thus, the degree of $a(x)$ is deg $a(x)=p-2$. The derivative of a constant no longer vanishes and the degree of the polynomial function is preserved! Over GF(7), for instance, we deal with polynomials of degree 5 and assume

$\frac{d1}{dx} \equiv \frac{d(-x^6)}{dx} \equiv x^5 \mod 7$ (negacyclic)

whereas

$\frac{d1}{dx} \equiv \frac{d(x^7)}{dx} \equiv 0 \mod 7$ (Hasse).

## 5. β-ADIC EXPANSIONS AND FINITE FIELD TAYLOR SERIES

The classical Taylor series of a function $a(x)$ is given by

$a(x) = a(x_0) + \frac{a^{(1)}(x_0)}{1!}(x-x_0) + \frac{a^{(2)}(x_0)}{2!}(x-x_0)^2 + \ldots$
$a^{(N-1)}(x_0)\frac{(x-x_o)^{N-1}}{(N-1)!} + R_n$

where $R_n = a^{(N)}(\zeta)\frac{(x-x_o)^N}{N!}$ ($x_0 < \zeta < x$)

(Lagrange's remainder).

On a finite field, Taylor series developed around an arbitrary point $\beta \in GF(p)$ can be considered according to

$a(x) = a(\beta) + a^{[1]}(\beta)(x-\beta) +$
$a^{[2]}(\beta)(x-\beta)^2 + \ldots + a^{[p-1]}(\beta)(x-\beta)^{p-1}$.

It is interesting that working with $x^p - 1 \equiv 0$ i.e., deg $a(x)=p-1$, we have

$R_p = \frac{a^{(p)}(\zeta)}{p!}(x-\beta)^p \equiv 0$ (Section 4) and the series is finite (no rounded error, although $\zeta$ is unknown).

Now let the β-adic expansion of $a(x)$ be
$a(x) \equiv b_0 + b_1(x-\beta) + b_2(x-\beta)^2 + \ldots + b_{N-1}(x-\beta)^{N-1}$
with $b_0 = a(\beta)$. It follows from Definition 2 that

$a^{[r]}(x) := \sum_{i=r}^{N-1} \binom{i}{r} b_i x^{i-r}$ so that $a^{[r]}(\beta) \equiv \binom{r}{r} b_r$

The unicity of the Taylor Finite Field series follows directly from [12]:

<u>Unicity</u>. *Let K be a field and K[x] be a ring of polynomials over such a field. The adic expansion of a polynomial a(x) over a ring K[x] is unique.*

Any function $a(x)$ over GF($p$) of degree $N$-1 can be written, for any element $\beta \in$ GF($p$), as
$$a(x) = \sum_{i=0}^{N-1} a^{[i]}(\beta)(x-\beta)^i \quad \text{without error or}$$
approximations. For instance, the $2^x$ GF(3)-function can be developed as $2^x \equiv 1 + 2x + 2x^2 \mod 3$ (McLaurin series).
Therefore,
$(2^x)^{[1]} \equiv 2 + x \; ; \; (2^x)^{[2]} \equiv \frac{1}{2!} \equiv 2 \mod 3.$

Letting $2^x \equiv b_0 + b_1(x-\beta) + b_2(x-\beta)^2 \mod 3$
then $2^x \equiv (b_0 - \beta b_1 + \beta^2 b_2) + (b_1 - 2\beta b_2)x + b_2 x^2$
with $b_0 = 2^\beta \mod 3$.
The $\beta$-adic expansions mod $p$=3 are:
$2^x \equiv 1 + 2x + 2x^2$           0-adic
$2^x \equiv 2 + 2(x-1)^2$           1-adic
$2^x \equiv 1 + (x-2) + 2(x-2)^2$   2-adic.
Indeed, $\forall \beta$,
$$\frac{a(x) - a(\beta)}{x-\beta} \equiv a^{[1]}(\beta) + \sum_{i=2}^{N-1} a^{[i]}(\beta)(x-\beta)^{i-1}, \text{ so one}$$
can, *by an abuse*, interpret the Hasse derivative as a "solution" of a "0/0-indetermination" over GF($p$):
$$a^{[1]}(\beta) \stackrel{?}{=} \left. \frac{a(x) - a(\beta)}{x-\beta} \right|_\beta .$$

## 6. APPLICATIONS: EXPONENTIAL AND TRIGONOMETRIC FUNCTIONS

Let us consider classical (real) series:
$$e^x = 1 + x + \frac{x^2}{2!} + \frac{x^3}{3!} + ...$$
$$\cos(x) = 1 - \frac{x^2}{2!} + \frac{x^4}{4!} - \frac{x^6}{6!} + ...$$
$$\sin(x) = x - \frac{x^3}{3!} + \frac{x^5}{5!} - \frac{x^7}{7!} + ...$$

As an example, we consider the ring GF($p$)[$x$] with $x^{p-1} + 1 \equiv 0$ so negacyclic Hasse derivatives can be used. Since deg $a(x) = p-2$ we truncate the series yielding
$(e^x)_7 \equiv 1 + x + 4x^2 + 6x^3 + 5x^4 + x^5 \pmod 7$,
$(\cos x)_7 \equiv 1 + 4x^2 + 5x^4 \pmod 7$,
$(\sin x)_7 \equiv x + 6x^3 + x^5 \pmod 7$.
The finite field Euler constant e and the exponential function can be interpreted as follows:
$(e^x)_p^{[1]} \equiv (e^x)_p$,
so that $(e^x)_p^{[r]} \equiv (e^x)_p$, $r \geq 1$.
$e \equiv 1 + \frac{1}{1!} + \frac{1}{2!} + \frac{1}{3!} + ... + \frac{1}{(p-2)!}$ (mod p).
Therefore $e^i \equiv 1 + i + 4i^2 + 6i^3 + 5i^4 + i^5 \mod 7$.
$i=$   0   1   2   3   4   5
$e^i =$   1   4   4   6   6   1.

In order to introduce trigonometric functions, it is possible to consider the complex $j = \sqrt{-1}$ since -1 is a quadratic non-residue for $p \equiv 3 \pmod 4$ and pick $e^{j \cdot i}$ so that new "e-trigonometric" sin and cos functions are:
$\cos(i) = \Re \exp(j.i) \equiv 1 + 4i^2 + 5i^4 \pmod 7$,
$\sin(i) = \Im \exp(j.i) \equiv i + 6i^3 + i^5 \pmod 7$.

It is straightforward to verify that $\cos(i)$ is even and $\sin(i)$ is an odd function. The negacyclic Hasse derivative yields:
$\sin^{[1]} i = \cos(i)$ and $\cos^{[1]} i = -\sin(i)$.

We have two kinds of trigonometric functions over a finite field: the $k$-trigonometric ones and these new e-trigonometric functions. Unfortunately, $\cos(i)$ and $\sin(i)$ do not lie on the unit circle, i.e., $\sin^2 i + \cos^2 i \not\equiv 1 \pmod p$. The "number e" is not a "natural" element of the field (in contrast with the real case). However, the $k$-trigonometric case is defined using an element of the field and presents a lot of properties. Nevertheless, derivatives formulae of $k$-trigonometric functions are not related as usual trigonometric derivatives.

These results can be generalized so as to define *hyperbolic* functions over finite fields, $p$ odd, according to the series developments:
$\forall i \in$ GF($p$)
$$\cosh(i) := 1 + \frac{i^2}{2!} + \frac{i^4}{4!} + \frac{i^6}{6!} + ... + \frac{i^{p-3}}{(p-3)!} \text{ and}$$
$$\sinh(i) := i + \frac{i^3}{3!} + \frac{i^5}{5!} + \frac{i^7}{7!} + ... + \frac{i^{p-2}}{(p-2)!}.$$
Indeed $\sinh(i) + \cosh(i) \equiv \exp(i)$.

We recall the convention $\alpha^{-\infty} = 0$ ($\forall \alpha) \in$ GF($p$) [7]. The idea is to consider an extended Galois Field by appending a symbol $-\infty$ as done for the Real number set ([1], def.1.39). ($\forall x) \in$ GF($p$):
$x + (-\infty) \equiv -\infty; \; x.(-\infty) \equiv -\infty$ and $x/(-\infty) \equiv 0$.

Therefore, other functions such as an hyperbole $1/x$ can also be defined.
$i=$   0   1   2   3   4
$1/i =$   $-\infty$   1   3   2   4.

Finally, it is interesting to consider "log" functions over a Galois field as the inverse of exponential. We have then
$i=$   0   1   2   3   4
$\log_2 i =$   $-\infty$   0   1   3   2.

$\log_2(0)$ and $\log_2(1)$ give the expected values. MacLaurin series cannot be derived as usual. Moreover, the $(e^x)$ function does not have an inverse so (ln$x$) is not defined, although $\log_\alpha x$ exists, for any primitive element $\alpha$.

# 7. INTEGRAL OVER A FINITE FIELD

We evaluate the sum of the $n^{th}$ power of all the $p$ elements of GF($p$). It is also assumed that $0^0 := 1$, since that $f(0)=a_0$.

$$0 + 1 + 2 + 3 + \ldots + (p-1) \equiv ?$$
$$0 + 1^2 + 2^2 + 3^2 + \ldots + (p-1)^2 \equiv ?$$
$$\ldots$$
$$0 + 1^{p-1} + 2^{p-1} + 3^{p-1} + \ldots + (p-1)^{p-1} \equiv ?$$

Lemma 1. ($n^{th}$ power summation).
$$\sum_{x=0}^{p-1} x^n \equiv \begin{cases} p-1 & \text{if } n = p-1 \\ 0 & \text{otherwise} \end{cases}.$$

proof. Letting $S_n := 1^n + 2^n + 3^n + \ldots + (p-1)^n$ $n=0,1,2,3\ldots,(p-1)$. Trivial case $n=0$ apart,
$$\binom{n+1}{1}S_1 + \binom{n+1}{2}S_2 + \ldots + \binom{n+1}{n}S_n = p^{n+1} - p \equiv 0$$
so that $n \neq p-1$ implies $S_1 \equiv S_2 \equiv \ldots \equiv S_n \equiv 0 \pmod{p}$. ∎

Definition 3: The definite integral over GF($p$) of a function $f$ from 0 to ($p$-1) is
$$I = \int_0^{p-1} f(\alpha)d\alpha := \sum_{\alpha=0}^{p-1} f(\alpha).$$

Proposition 2. The integral of a function $f$ presenting a series $a(x) = \sum_{i=0}^{p-1} a_i x^i$ is given by
$$\sum_{x=0}^{p-1} f(x) \equiv (p-1)a_{p-1} \mod p.$$

proof. Substituting the series of $f$ in the integral definition and applying Lemma 1,
$$I = \sum_{\alpha=0}^{p-1}\sum_{i=0}^{p-1} a_i \alpha^i = \sum_{i=0}^{p-1} a_i \sum_{x=0}^{p-1} x^i =$$
$$= \sum_{i=0}^{p-1} a_i (p-1)\delta_{i,p-1}$$
where $\delta_{i,j}$ is the Kronecker delta (1 if $i=j$, 0 otherwise). ∎

Corollary. If a function $f$ defined $\forall \alpha \in$ GF($p$) admits an inverse, then the $(p-1)^{th}$ coefficient of its series vanishes.

The condition $a_{p-1} \equiv 0 \mod p$ is a necessary condition but not a sufficient one to guarantee that $f$ is invertible. In fact the integral vanishes since all the images are distinct (the sum of all the elements of the field).

Proposition 3. Let $f$ and $g$ be functions defined for every element of GF($p$). Denoting by $\{a_i\}$ and $\{b_i\}$ the coefficients of their respective MacLaurin series, then

$$\int_0^{p-1} f \cdot g(\alpha)d\alpha = \sum_{k=0}^{p-1} f(k)g(k) \equiv \frac{1}{p-1}\sum_{i=0}^{p-1} a_i b_{p-1-i}$$

proof. Substituting the series of $f$ and $g$ in the sum and changing the order of sums, we have
$$\sum_{k=0}^{p-1} f(k)g(k) = \sum_{i,j=0}^{p-1} a_i b_j \sum_{k=0}^{p-1} k^{i+j}.$$
The proof follows from Lemma 1 and the fact that the inverse $(p-1)^{-1} \equiv p-1 \mod p$. ∎

In order to evaluate the integral over another "interval of GF($p$)", we consider:
$$\int_0^N f(\alpha)d\alpha := \sum_{\alpha=0}^N f(\alpha) = \sum_{i=0}^{p-1} a_i S_N(i),$$
where $S_N(i) := \sum_{x=0}^N x^i$.

A table with the values of $S_N(i)$ over GF(5) is showed below.

|          | $i=0$ | $i=1$ | $i=2$ | $i=3$ | $i=4$ |
|----------|-------|-------|-------|-------|-------|
| $S_0(i)$ | 1     | 0     | 0     | 0     | 0     |
| $S_1(i)$ | 2     | 1     | 1     | 1     | 1     |
| $S_2(i)$ | 3     | 3     | 0     | 4     | 2     |
| $S_3(i)$ | 4     | 1     | 4     | 1     | 3     |
| $S_4(i)$ | 0     | 0     | 0     | 0     | 4     |

Another interesting result establishes a link between Hasse derivatives and the Finite Field Fourier Transform [13].

The Galois field function obtained by an image element permutation of $f$ is here referred to as an $f$-permutation.

Proposition 4. The coefficients of the McLaurin series expansion of a given signal $f$ are exactly the inverse finite field Fourier transform (FFFT) of an $f$-permutation.

Proof. Given the 0-adic development $a(x)$ of a discrete function $f$, we consider then the finite field Fourier transform pair: $a \times A$ where the FFFT is from GF($p$) to GF($p$), and $A_i = a(\alpha^i) = \sum_{j=0}^{p-1} a_j \alpha^{ij}$. Since $a(x)$ is used to interpolate the discrete function $f(.)$, for any given $i$ there exist $j \in \{-\infty, 0, 1, \ldots, p-2\}$ such that $f(\alpha^j) = a(\alpha^i) = A_i$. ∎

# 8. CONCLUSIONS

The main point of this paper is to introduce the background of new tools useful for Engineering applications involving finite fields. A number of discrete functions is analyzed including a brief look at transcendental functions over Galois Fields. Finite field (FF) derivatives are considered, which in turn leads to Finite Field Taylor Series. A brief introduction to FF integration is also presented. Clearly many other extensions and applications do exist.

**Acknowledgments**- *The authors wish to thank Prof. V.C. da Rocha Jr. (Universidade Federal de Pernambuco) for introducing us Hasse derivatives.*